\documentclass[11pt]{article}
\usepackage{cite}
\usepackage{mathrsfs}
\usepackage{amssymb,amsfonts,amsmath}
\usepackage{enumerate}
\usepackage{indentfirst}
\usepackage{latexsym,bm}
\usepackage[noend]{algpseudocode}
\usepackage{algorithmicx,algorithm}
\usepackage{algorithm}
\usepackage{makeidx}
\usepackage{fancybox}
\usepackage{color}
\usepackage{multicol}
\usepackage[latin1]{inputenc}
\usepackage{graphicx}
\usepackage{epstopdf}
\usepackage{pstricks,pst-node,pst-text,pst-3d}
\usepackage{tikz}
\newtheorem{thm}{Theorem}[section]

\newtheorem{lem}[thm]{Lemma}

\newtheorem{obs}[thm]{Observation}

\newenvironment{pf}{{\noindent \it \bf Proof:}}{{\hfill$\Box$}\\}

\def\qed{\hfill \nopagebreak\rule{5pt}{8pt}}

\begin{document}

\title{\bf Kings and Kernels in Semicomplete Compositions}
\author{
\small Yuefang Sun$^{1}$, Zemin Jin$^{2}$\footnote{Corresponding author: Z.M. Jin
\vskip1mm Email address:
sunyuefang@nbu.edu.cn (Sun), zeminjin@zjnu.cn (Jin)} \\
\small $^{1}$ School of Mathematics and Statistics, Ningbo University,\\
\small Ningbo 315211, P. R. China \\
\small $^{2}$ Department of Mathematics, Zhejiang Normal University,\\
\small Jinhua 321004, P. R. China
}

\date{}
\maketitle

\begin{abstract}
Let $k$ be an integer with $k\geq 2$. A $k$-king in a digraph $D$ is a vertex which can reach every other vertex by a directed path of length at most $k$ and a non-king is a vertex which is not a 3-king.
A subset $K$ is $k$-independent if for every pair of vertices $x,y \in K$, we have $d_D(x, y), d_D(y, x)\geq k$; it is $\ell$-absorbent if for every $x\in V(D)\setminus K$ there exists $y\in K$ such that $d_D(x, y)\leq \ell$. A $k$-kernel of $D$ is a $k$-independent and $(k-1)$-absorbent subset of $V(D)$. A kernel is a 2-kernel. A set $K\subseteq V(D)$ is a quasi-kernel of $D$ if it is independent, and for every vertex $x\in V(D)\setminus K$, there exists $y\in K$ such that $d_D(x, y)\leq 2$. The problem {\sc $k$-Kernel} is determining whether a given digraph has a $k$-kernel. 

Let $Q=T[H_1, \dots, H_t]$ be the composition of $T$ and $H_i$ ($1\leq i\leq t, t\ge 2$), where $T$ is a digraph with $t$ vertices, and $H_1, \dots, H_t$ are pairwise disjoint digraphs. The composition $Q=T[H_1, \dots, H_t]$ is a semicomplete composition if $T$ is semicomplete. In this paper, we study kings and kernels in semicomplete compositions.

For the topic of kings, we characterize digraph compositions with a $k$-king and digraph compositions all of whose vertices are $k$-kings, respectively. We also discuss the existence of 3-kings, and study the minimum number of 4-kings in a strong semicomplete composition.

For the topic of kernels, we first 
study the existence of a pair of disjoint quasi-kernels in semicomplete compositions. We then deduce that the problem {\sc $k$-Kernel} restricted to strong semicomplete compositions is NP-complete when $k\in \{2,3\}$, and is polynomial-time solvable when $k\geq 4$. We also prove that when $k$ is divisible by 2 or 3, 
the problem {\sc $k$-Kernel} restricted to non-strong semicomplete compositions is NP-complete.
\vspace{0.3cm}
\vspace{0.3cm}

\end{abstract}

\section{Introduction}
We refer the readers to \cite{Bang-Jensen-Gutin, Bondy} for graph-theoretical notation and terminology not given here. All digraphs considered in this paper have no parallel arcs
or loops, and paths and cycles are always assumed to be directed, unless otherwise stated. We use $[n]$ to denote the set of all natural numbers from 1 to $n$.

Let $D$ be a digraph. If there is an arc from a vertex $x$ to a vertex $y$ in $D$, then $x$ {\em dominates} $y$, denoted by $x\rightarrow y$. We shall use the notation $x\not\rightarrow y$ if there is no arc from $x$ to $y$. If $A$ and $B$ are two subdigraphs of $D$ and every vertex of $A$ dominates each vertex of $B$, then $A$ {\em dominates} $B$, denoted by $A\rightarrow B$.
We shall use $A\Rightarrow B$ to denote that $A$ dominates $B$ and there is no arc from $B$ to $A$. A {\em strong component} of a digraph $D$ is a maximal induced subdigraph
of $D$ which is strong. 
An {\em initial strong component} is one which has no arcs coming in from any other strong components.
A vertex $x$ in $D$ is a {\em sink} (resp. {\em source}) if $d^+(x)=0$ (resp. $d^-(x)=0$). For every pair of distinct vertices $x, y\in V(D)$, we use $d_D(x, y)$ to denote the length of a shortest path from $x$ to $y$ (if there is no such path, then $d_D(x, y)$ is infinite). For a vertex $x$ and a vertex subset $V'$ of $V(D)$, let $d_D(x, V')=\min\{d_D(x,y)\colon\ y\in V'\}$.

A digraph $D$ is {\em vertex-pancyclic} if for every $x\in V(D)$ and every integer $k$ with $3\leq k\leq n$, there exists a $k$-cycle through $x$ in $D$.
A {\em transitive tournament} with vertex set $\{u_i\colon\ i\in [t]\}$ is a tournament such that  $u_i\rightarrow u_j$ but $u_j\not\rightarrow u_i$ for any pair $i, j$ with $i, j\in [t]$. A digraph is {\em semicomplete multipartite} if it is obtained from a complete multipartite graph by replacing every edge by an arc or a pair of opposite arcs. A {\em multipartite~tournament} is a semicomplete multipartite digraph with no 2-cycle.
Let $T$ be a digraph with vertices $u_1, \dots, u_t$ ($t\ge 2$), and let $H_1, \dots, H_t$ be digraphs which are pairwise vertex disjoint. 
The {\em composition} $Q=T[H_1, \dots, H_t]$ is a digraph with vertex set $V(H_1)\cup V(H_2)\cup \ldots \cup V(H_t)$ and arc set $(\cup_{i=1}^t{A(H_i)})\cup \{h_ih_j\colon\ h_i\in V(H_i), h_j\in V(H_j), u_iu_j\in A(T)\}$. 
The composition $Q=T[H_1, \dots, H_t]$ is a {\em semicomplete composition} if $T$ is semicomplete. If $Q=T[H_1, \dots, H_t]$ and none of $H_1, \dots, H_t$ has an arc, then $Q$ is an {\em extension} of $T$. For any set $\Phi$ of digraphs, $\Phi^{ext}$ denotes the (infinite) set of all extensions of digraphs in $\Phi$, which are {\em extended $\Phi$-digraphs}. By definition, the class of extended semicomplete digraphs is a subclass of semicomplete compositions. 
A digraph $D$ is {\em quasi-transitive}, if for any triple $x,y,z$ of distinct vertices of $D$, if $xy$ and $yz$ are arcs of $D$
then either $xz$ or $zx$ or both are arcs of $D.$
Bang-Jensen and Huang gave a complete characterization of quasi-transitive digraphs by using digraph compositions as shown below.

\begin{thm}\label{intro01}\cite{Bang-Jensen-Huang1}
Let $D$ be a quasi-transitive digraph. The following assertions hold:
\begin{description}
\item[(a)] If $D$ is not strong, then there exists a transitive oriented
graph $T$ with vertices $\{u_i\colon\ i\in [t]\}$ and strong quasi-transitive digraphs $H_1, \dots, H_t$ such that $D = T[H_1, \dots, H_t]$, where $H_i$ is substituted for $u_i$, $i\in [t]$.
\item[(b)] If $D$ is strong, then there exists a strong semicomplete
digraph $S$ with vertices $\{v_j\colon\ j\in [s]\}$ and quasi-transitive digraphs $Q_1, \dots, Q_s$ such that $Q_j$ is either a vertex or is non-strong and $D = S[Q_1, \dots, Q_s]$, where $Q_j$ is substituted for $v_j$, $j\in [s]$.
\end{description}
\end{thm}

Digraph compositions generalize some families of digraphs, including quasi-transitive digraphs (by Theorem~\ref{intro01}), extended semicomplete digraphs and lexicographic product digraphs (when $H_i$ is the same digraph $H$ for every $i\in [t]$, $Q$ is the lexicographic product of $T$ and $H$, see, e.g., \cite{Hammack}). In particular, strong semicomplete compositions generalize strong quasi-transitive digraphs (by Theorem~\ref{intro01}). To see that strong semicomplete compositions form a significant generalization of strong quasi-transitive digraphs, observe that the Hamiltonian cycle problem is polynomial-time solvable for quasi-transitive digraphs \cite{Gutin4}, but NP-complete for strong semicomplete compositions (see, e.g., \cite{Bang-Jensen-Gutin-Yeo}).
While digraph composition has been used since 1990s to study quasi-transitive digraphs and their generalizations, see, e.g., \cite{Bang-Jensen-Gutin, Bang-Jensen-Huang1}, the study of digraph decompositions in their own right was initiated  by Sun, Gutin and Ai in \cite{Sun-Gutin-Ai}, where some results on the existence of arc-disjoint strong spanning subdigraphs were given. After that, Bang-Jensen, Gutin and Yeo \cite{Bang-Jensen-Gutin-Yeo} characterized all strong semicomplete compositions with a pair of arc-disjoint strong spanning subdigraphs. Gutin and Sun \cite{Gutin-Sun} studied the existence of a pair of arc-disjoint in- and out-branchings rooted at the same vertex in this class of digraphs. Especially, they characterized semicomplete compositions with such a pair rooted at the same vertex, which generalizes the corresponding characterization by Bang-Jensen and Huang \cite{Bang-Jensen-Huang1} for quasi-transitive digraphs. In \cite{Sun-Gutin-Zhang}, Sun, Gutin and Zhang studied the strong subgraph packing problem in compositions of digraphs. Recently, Sun and Jin \cite{Sun-Jin} studied the structural properties of semicomplete compositions, including results on connectivity, paths, cycles, strong spanning subdigraphs and acyclic spanning subdigraphs. These results show that this class of digraphs shares some nice properties of other classes of digraphs, like quasi-transitive digraphs.

A {\em $k$-king} in a digraph $D$ is a vertex which can reach every other vertex by a directed path of length at most $k$. A $k$-king $u$ is {\em strict} if there exists a vertex $v$ such that $d_D(u,v)=k$. A 2-king is a {\em king}, and a {\em non-king} is a vertex which is not a 3-king. The study of kings in digraphs began with Landau \cite{Landau} and the term king was introduced by Maurer \cite{Maurer}. The concept of $k$-king was first introduced in \cite{Bang-Jensen-Huang2} when the authors studied quasi-transitive digraphs, and it has attracted much attention from researchers \cite{Bang-Jensen-Gutin}. 

Let $D$ be a digraph. A set $K\subseteq V(D)$ is a {\em kernel} of $D$ if it is an independent set such that every vertex
in $V(D)\setminus K$ dominates some vertex in $K$. von Neumann and Morgenstern introduced the concept of kernel when they studied cooperative games \cite{Neumann-Morgenstern}, and after that, researchers from other areas, such as list coloring, game theory and perfect graphs \cite{Boros-Gurvich}, mathematical logic \cite{Bezem-Grabmayer-Walicki} and complexity theory \cite{Walicki-Dyrkolbotn}, studied this topic.
There are many generalizations of this concept. A set $K\subseteq V(D)$ is a {\em quasi-kernel} of $D$ if it is independent, and for every vertex $x\in V(D)\setminus K$, there exists $y\in K$ such that $d_D(x, y)\leq 2$. A subset $K$ is {\em $k$-independent} if for every pair of vertices $x,y \in K$, we have $d_D(x, y), d_D(y, x)\geq k$; it is {\em $\ell$-absorbent} if for every $x\in V(D)\setminus K$ there exists $y\in K$ such that $d_D(x, y)\leq \ell$. A {\em $k$-kernel} of $D$ is a $k$-independent and $(k-1)$-absorbent subset of $V(D)$. Clearly, a 2-kernel is exactly a kernel. The problem {\sc $k$-Kernel} is determining whether a given digraph has a $k$-kernel.

In this paper, we continue the research of digraph compositions and focus on the topics of kings and kernels in semicomplete compositions. In Section~\ref{sec:king}, the existence of $k$-kings are discussed. We first give characterizations of all digraph compositions with a $k$-king and digraph compositions all of whose vertices are $k$-kings, respectively, where $k\geq 2$ is an integer (Theorem~\ref{thma}). When $T$ is semicomplete, we first give a sufficient condition to guarantee that $Q$ has at least two 3-kings. In particular,  every strong semicomplete composition has at least two 3-kings (Theorem~\ref{thmd}). The adjacency between 3-kings and non-kings in a strong semicomplete composition is also studied (Lemma~\ref{thmc}). In this paper, for a semicomplete composition $Q$, we give a sufficient condition under which there exists a semicomplete composition $Q'$ which contains $Q$ as an induced subdigraph such that the set of all 3-kings of $Q'$ is precisely $V(Q)$ (Theorem~\ref{thme}). Similar problems for tournaments and quasi-transitive digraphs have been studied in \cite{Bang-Jensen-Huang2, Huang-Li, Reid}. In the end of this section, we discuss the minimum number of 4-kings in a strong semicomplete composition (Theorem~\ref{thmf}).


The topic of kernels is discussed in Section~\ref{sec:Kernel}. 
Gutin et al. \cite{Gutin-Koh-Tay-Yeo} characterized digraphs with exactly one or two quasi-kernels. In particular, if a digraph has precisely two quasi-kernels then these two quasi-kernels are actually disjoint. This raises the question of which digraphs contain a pair of disjoint quasi-kernels. We deduce that when $T$ has no sink, each semicomplete composition contains a pair of disjoint quasi-kernels. In particular, every strong semicomplete composition contains a pair of disjoint quasi-kernels (Theorem~\ref{thmb1}). It is known that the problems {\sc $2$-Kernel} (also called {\sc Kernel}) and {\sc $3$-Kernel} are both NP-Complete (see e.g. \cite{Bang-Jensen-Gutin, Bang-Jensen-Gutin2}). In the end of Section~\ref{sec:Kernel}, we study the complexity of {\sc $k$-Kernel} for any integer $k\geq 2$ and deduce that this problem restricted to strong semicomplete compositions is NP-complete when $k\in \{2,3\}$ and is polynomial-time solvable when $k\geq 4$ (Theorem~\ref{thmd2}(a)). We also prove that when $k$ is divisible by 2 or 3, the problem {\sc $k$-Kernel} on general digraphs is NP-complete (Lemma~\ref{lem-kkernel}), by this result, we prove that the problem {\sc $k$-Kernel} restricted to non-strong semicomplete compositions is NP-complete when $k$ is divisible by 2 or 3 (Theorem~\ref{thmd2}(b)).

\section{Kings}\label{sec:king}

\subsection{$T$ is arbitrary}

By definition, we directly have the following observations.

\begin{obs}\label{obs3}
Let $Q=T[H_1, \dots, H_t]$ be a digraph composition, $Q$ is strong if and only if $T$ is strong.
\end{obs}

\begin{obs}\label{obs2}
Let $k\geq 2$ be an integer. If $H$ is a spanning subdigraph of $D$, then every $k$-king (if exists) of $H$ is also a $k$-king of $D$.
\end{obs}

\begin{obs}\label{obs1}
Let $k\geq 2$ be an integer. If $Q=T[H_1, \dots, H_t]$ has a $k$-king $v\in V(H_i)$ for some $i \in [t]$, then $u_i$ is a $k$-king of $T$.
\end{obs}

For any integer $k\geq 2$, we characterize all digraph compositions with a $k$-king and digraph compositions all of whose vertices are $k$-kings, respectively.

\begin{thm}\label{thma}
Let $k\geq 2$ be an integer, and $Q=T[H_1, \dots, H_t]$ be a digraph composition. The following assertions hold:
\begin{description}
\item[(a)] $Q$ has a $k$-king if and only if $T$ has a $k$-king $u_{i}$ for some $i\in [t]$ and at least one of the following possibilities holds: $(i)$~$H_{i}$ has a $k$-king; $(ii)$~$|V(H_{i})|\geq 2$ and $u_{i}$ belongs to a cycle of $T$ with length at most $k$.
\item[(b)] All vertices of $Q$ are $k$-kings if and only if for each $i\in [t]$, $u_{i}$ is a $k$-king of $T$ and at least one of the following possibilities holds: $(i)$~all vertices of $H_{i}$ are $k$-kings of $H_i$;   $(ii)$~$|V(H_{i})|\geq 2$ and $u_{i}$ belongs to a cycle of $T$ with length at most $k$.
\end{description}
\end{thm}
\begin{pf}

\noindent{\bf Part (a)} The proof of {\em sufficiency} is not hard, so we just prove the {\em necessity}. If $Q=T[H_1, \dots, H_t]$ has a $k$-king $v$ which belongs to some $V(H_i)$, then $u_i$ is a $k$-king of $T$ by Observation \ref{obs1}. If $|V(H_i)|=1$, then the unique vertex $v$ in $H_i$ is of course a $k$-king of $H_i$. Hence in the following we assume that $|V(H_i)|\geq 2$.

If $H_i$ has no $k$-king, then since $v$ is a $k$-king of $Q$, there must exist a vertex $x\in V(H_i)\setminus \{v\}$ such that there is a $v,x$-path, denoted by $P_{vx}$, of length at most $k$ which uses at least one vertex from $V(Q)\setminus V(H_i)$. We claim that $u_i$ belongs to a cycle of length at most $k$ in $T$. To see this, we replace the vertices in $P_{vx}$ (in order) with their corresponding vertices in $T$. For example, we replace $v$ with $u_i$. We can get a closed walk of length at most $k$ which contains $u_i$ in $T$, this argument means that $u_i$ indeed belongs to a cycle of length at most $k$ in $T$. Therefore, either $H_{i}$ has a $k$-king, or, $|V(H_{i})|\geq 2$ and $u_{i}$ belongs to a cycle of $T$ with length at most $k$.

\noindent{\bf Part (b)} The proof of {\em sufficiency} is not hard and we just prove the {\em necessity}. By Observation \ref{obs1}, all vertices of $T$ are $k$-kings of $T$. For any $i\in [t]$, if $|V(H_i)|=1$, then of course the unique vertex is a $k$-king of $H_i$. In the following we assume that $|V(H_i)|\geq 2$. If one vertex, say $u$, of $H_i$ is not $k$-king of $H_i$, then (since $u$ is a $k$-king of $Q$) there must exist a vertex $x\in V(H_i)\setminus \{u\}$ such that there is a $u,x$-path, denoted by $P_{ux}$, of length at most $k$ which uses at least one vertex from $V(Q)\setminus V(H_i)$. With a similar argument to that of \noindent{\bf Part (a)}, we deduce that $u_i$ belongs to a cycle of length at most $k$ in $T$.
\end{pf}



\subsection{$T$ is semicomplete}



The following result can be found in literature, see e.g. Theorem~2.2.9 of \cite{Bang-Jensen-Gutin2}.

\begin{thm}\label{thm02}
Every strong semicomplete digraph is vertex-pancyclic.
\end{thm}

Landau got the following result which is used in our following results, such as Theorem~\ref{thmd}.

\begin{thm}\label{thm01}\cite{Landau}
Every tournament has a king. More precisely, every vertex with maximum out-degree is a king.
\end{thm}

Bang-Jensen and Huang obtained the following result on the existence of 3-kings in a quasi-transitive digraph.

\begin{thm}\label{Bang-Huangking}\cite{Bang-Jensen-Huang2}
Let $D$ be a quasi-transitive digraph with a 3-king. If $D$ has no source, then $D$ has at least two 3-kings.
\end{thm}



For semicomplete compositions, we give a similar result to Theorem~\ref{Bang-Huangking}.

\begin{thm}\label{thmd}
Let $Q=T[H_1, \dots, H_t]$ be a semicomplete composition with a 3-king. If $T$ has no source, then $Q$ has at least two 3-kings. In particular, every strong semicomplete composition $Q=T[H_1, \dots, H_t]$ has at least two 3-kings, and a vertex $v\in H_i$ is a 3-king of $Q$ if and only if $u_i$ is a 3-king of $T$. 
\end{thm}
\begin{pf} We first prove the following claim.

{\bf Claim.} Every strong semicomplete composition $Q=T[H_1, \dots, H_t]$ has a 3-king. Furthermore, a vertex $v\in H_i$ is a 3-king of $Q$ if and only if $u_i$ is a 3-king of $T$.

{\bf Proof of the claim.} By Theorem~\ref{thm01}, there is a 3-king, say $u_i$, in $T$. If $|V(H_i)|=1$, then the unique vertex in $V(H_i)$ is also a 3-king of $Q$. Hence we assume that $|V(H_i)|\geq 2$. By Theorem~\ref{thm02}, $u_i$ belongs to a 3-cycle in $T$, so $Q$ has a 3-king by Theorem \ref{thma}.

If $u_i$ is a 3-king of $T$, then $v\in H_i$ can reach every vertex outside $H_i$ by a path of length at most three. By Theorem~\ref{thm02}, $u_i$ belongs to a 3-cycle in $T$, so $v$ can reach the remaining vertices of $H_i$ by a path of length at most three. The other direction follows from Observation \ref{obs1}. \qed

\noindent{\bf Case 1:} $Q$ is strong.

By Observation~\ref{obs3}, $T$ is strong and of course has no source. By the above claim, $T$ has a 3-king and furthermore has at least two 3-kings by Theorem~\ref{Bang-Huangking}, hence $Q$ has at least two 3-kings, the conclusion holds by Claim.

\noindent{\bf Case 2:} $Q$ is not strong.

Since $Q$ has a 3-king, the initial strong component, say $Q'$, must be unique and contain all 3-kings. 
Without loss of generality, assume that $Q'$ contains vertices from $H_i~(i\in [t'], t'\geq 2)$. The component $Q'$ is a strong semicomplete composition with the form $Q'=T'[H_1, \dots, H_{t'}]$, where $V(T')=\{u_i\colon\ i\in [t']\}$ and $|T'|=t'\geq 2$. By Case 1, $Q'$ contains two 3-kings, say $x$ and $y$. Since $Q'$ is the unique initial strong component of $Q$, by the definition of semicomplete compositions, $x\rightarrow z$ and $y\rightarrow z$ for each $z\in V(Q)\setminus V(Q')$. Hence both $x$ and $y$ are 3-kings of $Q$.
\end{pf}

\noindent{\bf Remark A}: According to the second paragraph of the proof for the claim in Theorem~\ref{thmd}, we actually can divide $H_1, \dots, H_t$ in a strong semicomplete composition $Q=T[H_1, \dots, H_t]$ into those for which all $v\in H_i$ are 3-kings (in this case, $u_i$ is a 3-king of $T$) and those for which no $v\in H_i$ is a 3-king.

Recall that in Theorem~\ref{Bang-Huangking}, Bang-Jensen and Huang \cite{Bang-Jensen-Huang2} proved that a quasi-transitive digraph $D$ with a 3-king contains at least two 3-kings provided that it has no source.
It is worth noting that we cannot replace the condition ``If $T$ has no source'' in Theorem~\ref{thmd} by ``If $Q$ has no source'', according to the following example: 
let $Q=T[H_1, \dots, H_t]$ be a semicomplete composition such that $T$ is a transitive tournament with source $u_1$, and $H_1$ is a symmetric digraph with vertex set $\{v_{j}\colon\ 1\leq j\leq 7\}$ and arc set $\{v_{j}v_{j+1}, v_{j+1}v_{j} \colon\ 1\leq j\leq 6\}$. It can be checked that $Q$ has no source but $v_{4}$ is the unique 3-king of $Q$.

The following result concerns the adjacency between 3-kings and non-kings.

\begin{lem}\label{thmc}
Let $Q=T[H_1, \dots, H_t]$ be a strong semicomplete composition. There is an arc between every 3-king and every non-king; moreover, for every non-king $u$, there exists a 3-king $v$ such that $d_Q(u, v)>3$ and $v$ dominates $u$.
\end{lem}
\begin{pf} By Theorem~\ref{thmd} and Remark A, there is an arc between every 3-king and every non-king.
Let $u$ be any non-king of $Q$. Let $H$ be the subdigraph of $Q$ induced by the set of all vertices $v$ such that $d_Q(u, v)>3$. Notice that by Theorem~\ref{thm02}, $T$ is vertex-pancyclic and so $u$ can reach every other vertex of the same $H_i$ by a path of length at most three. We may assume without loss of generality that there is an integer $t'$ $(t'\in [t])$ such that $V(H)=\bigcup_{i=1}^{t'}{V(H_{i})}$. 

Without loss of generality, assume that $u_1$ is a king (by Theorem~\ref{thm01}) of the subdigraph $T'$ of $T$ induced by $\{u_i\colon\ i\in [t']\}$.
Let $x$ be a vertex of $H_1$. Suppose that $d_Q(x, y)>3$ for some vertex $y\in V(Q)$. By Theorem~\ref{thm02} and the fact that $u_1$ is a king of $T'$, we must have $y\not\in H$, and so $d_Q(u, y)\leq 3$.
Without loss of generality, assume that $d_Q(u, y)= 3$, and let $u, a, b, y$ be a shortest $u, y$-path. Observe that $a\not\in V(H_1)$, and hence there is an arc between $x$ and $a$. If $xa\in A(Q)$, then $x, a, b, y$ is an $x,y$-path of length 3, a contradiction. Otherwise, $ax\in A(Q)$, then $u, a, x$ is a $u,x$-path of length 2, which also produces a contradiction. Hence $x$ is a 3-king of $Q$ such that $xu\in A(Q)$ and $d_Q(u,x)>3$, as desired.
\end{pf}

Bang-Jensen and Huang \cite{Bang-Jensen-Huang2} proved that Lemma~\ref{thmc} holds for a quasi-transitive digraph with a 3-king, which means that the result may hold even for a non-strong quasi-transitive digraph (as long as it has a 3-king). However, for a non-strong semicomplete composition $Q$ with a 3-king, Lemma~\ref{thmc} may not hold. We just use the example before Lemma~\ref{thmc}: 
$Q$ is not strong and $v_{4}$ is the unique 3-king of $Q$; furthermore, there is no arc between $v_{4}$ and $v_{1}$, and of course $v_{4}$ does not dominate $v_{1}$.

For a semicomplete composition $Q$, we give a sufficient condition under which there exists a semicomplete composition $Q'$ which contains $Q$ as an induced subdigraph such that the set of all 3-kings of $Q'$ is precisely $V(Q)$. It suffices to study strong semicomplete compositions with a non-king as the case that all vertices are 3-kings is trivial.

\begin{thm}\label{thme}
Let $Q=T[H_1, \dots, H_t]$ be a strong semicomplete composition. If $T$ has a strict 3-king and every 2-king of $T$ is dominated by some strict 3-king of $T$, then there exists a semicomplete composition $Q'$ which contains $Q$ as an induced subdigraph such that the set of all 3-kings of $Q'$ is precisely $V(Q)$.
\end{thm}
\begin{pf}
Without loss of generality, let $\{u_i\colon\ i\in [p]\}$ and $\{u_i\colon\ p+1\leq i\leq q\}$ be the sets of strict 3-kings and 2-kings of $T$, respectively. Let $T'$ be a new semicomplete digraph satisfying: $(i)$ $V(T')=V(T)\cup \{v_i\colon\ i\in [p]\}$, $(ii)$ for each $i\in [p]$, $v_i\rightarrow u_i$ and $u\rightarrow v_i$ for any $u\in V(T)\setminus \{u_i\}$, $(iii)$ for any pair $i, j$ with $i, j\in [p]$, $v_i\rightarrow v_j$ if and only if $u_i\rightarrow u_j$. Let $Q'=T'[H_1, \dots, H_t, v_1, \dots, v_p]$. Observe that $Q'$ is strong.



Our result holds by the following two claims.

{\bf Claim 1.} For each $i\in [p]$, $v_i$ is a non-king of $Q'$.

{\bf Proof of Claim 1.} For each $i\in [p]$, observe that all vertices in $H_i$ are strict 3-kings of $Q$, so they cannot reach the vertices of some $H_j$ in two steps, hence $v_i$ cannot reach the vertices of $H_j$ in three steps in $Q'$ and so $v_i$ is a non-king of $Q'$.\qed

{\bf Claim 2.} All vertices of $Q$ are 3-kings of $Q'$.

{\bf Proof of Claim 2.} Let $u\in V(Q)$. We need to consider the following three cases:

The first case is that $u\in \bigcup_{i=1}^p{V(H_i)}$. Without loss of generality, assume that $u\in V(H_1)$. Since $u_{1}$ is a 3-king of $T$, $u$ can reach any vertex $v\in V(Q)\setminus V(H_1)$ in at most three steps. For any vertex $v\in V(H_1)\setminus \{u\}$, the path $u, w, v_1, v$ is a $u,v$-path of length three, where $w \in V(Q)\setminus V(H_1)$ (note that such a vertex $w$ exists since $u_{1}$ is a 3-king of $T$ and hence has an out-neighbour in $V(T)\setminus \{u_1\}$). Moreover, $u$ dominates all new vertices but $v_1$.

The second case is that $u\in \bigcup_{i=p+1}^q{V(H_i)}$. Without loss of generality, assume that $u\in V(H_{p+1})$. Since $u_{p+1}$ is a 2-king of $T$, $u$ can reach any vertex $v\in V(Q)\setminus V(H_{p+1})$ in at most two steps. 
By the assumption that every 2-king of $T$ is dominated by some strict 3-king of $T$, $u_{p+1}$ is dominated by some vertex, say $u_1$, from $\{u_i\colon\ i\in [p]\}$. Hence for any vertex $v\in V(H_{p+1})\setminus \{u\}$, the path $u, v_1, w, v$ is a $u,v$-path of length three, where $w\in V(H_1)$. Moreover, $u$ dominates all new vertices.

Finally we consider the case that $u$ is a non-king of $Q$. If $v\in \bigcup_{i=1}^p{V(H_i)}$, then $u$ can reach $v$ by a path of the form $u, v_j, v$ if $v\in V(H_j)$. If $v\in \bigcup_{i=p+1}^q{V(H_i)}$,
then $u$ can reach $v$ by a path of the form $u, v_j, w, v$, where $w\in V(H_j)$ for some $j\in [p]$. Such a path exists by the assumption that every 2-king of $T$ is dominated by some strict 3-king of $T$. If $v$ is also a non-king, then $u$ can reach $v$ by a path of the form $u, v_j, w, v$, where $w\in V(H_j)$ for some $j\in [p]$. To see why such a path exists, it suffices to notice that every non-king of $T$ must be dominated by some strict 3-king of $T$: Otherwise, there must exist some non-king $x$ which dominates all strict 3-kings of $T$. Consider any vertex $y\in V(T)$, if $y$ is a 2-king of $T$, then there exists an $x,y$-path in $T$: $x, x', y$, where $x'$ is a strict 3-king and its existence follows from the assumption that every 2-king of $T$ is dominated by some strict 3-king of $T$; if $y$ is also a non-king, then there exists an $x,y$-path in $T$: $x, x', y$, where $x'$ is a 3-king and its existence follows from Lemma~\ref{thmc} as $T$ itself is a strong semicomplete composition. 
Hence $x$ is a 3-king of $T$, which produces a contradiction.

Thus, in each case $u$ is a 3-king of $Q'$, so the claim holds.
\end{pf}

At the end of this section, we discuss 4-kings in semicomplete compositions. In our argument, some results on semicomplete multipartite digraphs are needed. 

For strong semicomplete multipartite digraphs, Gutin and Yeo obtained the following result on the number of 4-kings.

\begin{thm}\label{thm03}\cite{Gutin-Yeo}
Every strong semicomplete multipartite digraph with at least six vertices has at least five 4-kings.
\end{thm}

The following result on semicomplete bipartite digraphs was first stated by Wang and Zhang \cite{Wang-Zhang}.

\begin{thm}\label{thm06}\cite{Wang-Zhang}
Let $D$ be a semicomplete bipartite digraph with a unique initial strong component. If there is no 3-king in $D$, then there are at least eight 4-kings in $D$.
\end{thm}

For multipartite tournaments, Koh and Tan got the following result.

\begin{thm}\label{thm04}\cite{Koh-Tan, Koh-Tan2}
If a multipartite tournament has a unique initial strong component and no 3-king, then it has at least eight 4-kings.
\end{thm}

It is worth noting that many results for strong multipartite tournaments also hold for strong semicomplete multipartite digraphs, due to the following result by Volkmann.

\begin{thm}\label{thm05}\cite{Volkmann}
Every strong semicomplete $c$-partite digraph with $c \geq 3$ contains a spanning strong oriented subdigraph.
\end{thm}

Now we are ready to deal with the minimum number of $4$-kings in a strong semicomplete composition in the following result.

\begin{thm}\label{thmf}
Every strong semicomplete composition $Q=T[H_1, \dots, H_t]$ with at least six vertices has at least five 4-kings. Furthermore, if $Q$ has no 3-king, then it has at least eight 4-kings.
\end{thm}
\begin{pf} Let $Q'$ be the spanning subdigraph of $Q$ by deleting all arcs inside each $H_i$. 
By Observation~\ref{obs3}, $Q'$ is a strong semicomplete $t$-partite digraph with at least six vertices and so has at least five 4-kings by Theorem~\ref{thm03}. Hence $Q$ has at least five 4-kings.

If $Q$ has no 3-king, then $Q'$ also has no 3-king by Observation~\ref{obs2}. In the following we will show that $Q'$ has at least eight 4-kings, and hence the result directly holds.
The case that $t=2$ holds by Theorem~\ref{thm06}. For the case that $t\geq 3$, by Theorem~\ref{thm05}, $Q'$ has a strong spanning subdigraph $Q''$ which is a strong $t$-partite tournament (without a 3-king), hence $Q''$ (and so $Q'$) has at least eight 4-kings by Theorem~\ref{thm04}.
\end{pf}

\section{Kernels}\label{sec:Kernel}

Chv\'{a}tal and Lov\'{a}sz obtained the following result on quasi-kernels.

\begin{thm}\label{thm012}\cite{Chvatal-Lovasz}
Every digraph contains a quasi-kernel.
\end{thm}

Heard and Huang \cite{Heard-Huang} gave sufficient conditions to guarantee the existence of a pair of disjoint quasi-kernels in several classes of digraphs, including semicomplete digraphs.

\begin{thm}\cite{Heard-Huang}\label{thm011}
Every semicomplete digraph $D$ with no sink contains two vertices
$x, y$ such that $\{x\}$ and $\{y\}$ are both quasi-kernels of $D$.
\end{thm}

We now get a sufficient condition to guarantee the existence of a pair of disjoint quasi-kernels in semicomplete compositions.

\begin{thm}\label{thmb1}
Let $Q=T[H_1, \dots, H_t]$ be a semicomplete composition such that $T$ has no sink. The digraph $Q$ contains a pair of disjoint quasi-kernels. In particular, every strong semicomplete composition contains a pair of disjoint quasi-kernels.
\end{thm}
\begin{pf}
Let $Q=T[H_1, \dots, H_t]$ be a semicomplete composition such that $T$ has no sink. 
By Theorem~\ref{thm011}, there are two vertices, say $u_i$ and $u_j$, such that $\{u_i\}$ and $\{u_j\}$ are both quasi-kernels of $T$. Let $Q_i$ and $Q_j$ be the quasi-kernels of $H_i$ and $H_j$, respectively (note that this can be guaranteed by Theorem~\ref{thm012}). It can be checked that $Q_i$ and $Q_j$ are disjoint quasi-kernels of $Q$.

If $Q=T[H_1, \dots, H_t]$ is a strong semicomplete composition, then $T$ is strong by Observation~\ref{obs3} and hence contains no sink. By the above argument, $Q$  contains a pair of disjoint quasi-kernels.
\end{pf}

Note that when $Q=T[H_1, \dots, H_t]$ is not strong, it may not contain a pair of disjoint quasi-kernels. We use the example in Section~3 of \cite{Gutin-Koh-Tay-Yeo}, let $Q=T[H_1, \dots, H_t]$, where each $H_i$ is a directed path of order two and $T$ is a tournament satisfying: for every pair vertices $u, v$ of $V(T)$, there exists a vertex $w$ such that $u\rightarrow w$ and $v\rightarrow w$. It was shown in \cite{Gutin-Koh-Tay-Yeo} that $Q$ does not contain a pair of disjoint quasi-kernels.

In the remaining of this section, we turn our attention to the complexity of the problem {\sc $k$-Kernel} restricted to semicomplete compositions. We start with the following result which provides a necessary condition for a semicomplete composition to have a $k$-kernel when $k\geq 3$.

\begin{lem}\label{thmc1}
Let $Q=T[H_1, \dots, H_t]$ be a semicomplete composition. If $Q$ contains a $k$-kernel with $k\geq 3$, then there is a vertex $v\in V(Q)$ such that $\{v\}$ is a $(k-1)$-absorbent set of $Q-(V(H_i)\setminus \{v\})$, where $v\in V(H_i)$ for some $i\in [t]$.
\end{lem}
\begin{pf}
Let $K$ be a $k$-kernel of $Q$ and $v$ be any vertex of $K$. Observe that $K\subseteq V(H_i)$ for some $i\in [t]$. Let $Q'=Q-(V(H_i)\setminus \{v\})$ and $u$ be any vertex of $V(Q)\setminus V(H_i)$. Since $Q$ is a semicomplete composition, $u$ and $v$ must be adjacent in $Q$. If $u \rightarrow v$, then we are done. In the following we assume that $v \rightarrow u$.
By the definition of $K$, there must be $x\in K$ such that $d_Q(u, x)\leq k-1$. If $v=x$, then we are done. In the following we assume that $v\neq x$.

We claim that $d_Q(u, x)=k-1$. Indeed, if there is a $u,x$-path, denoted by $P_{u, x}$, in $Q$ of length at most $k-2$, then the path $vP_{u, x}$ is a $v,x$-path in $Q$ of length at most $k-1$, which contradicts the $k$-independence of $K$. Hence there is a $u,x$-path of length $k-1$ in $Q$: $u_0=u, u_1, \dots, u_{k-1}=x$. If there exists some $j\in [k-2]$ such that $u_j\in V(H_i)$ (we further assume that $j$ is the smallest integer satisfying this property), then the path $u_0, u_1, \dots, u_{j-1}, v$ is a $u,v$-path in $Q'$ of length at most $k-1$. Otherwise, the path $u_0, u_1, \dots, u_{k-2}, v$ is a $u,v$-path in $Q'$ of length $k-1$. Therefore, in both cases, we conclude that $d_{Q'}(u, v)\leq k-1$. This completes the proof.
\end{pf}

Furthermore, when $k\geq 4$, we can characterize strong semicomplete compositions with a $k$-kernel.

\begin{lem}\label{thmd1}
Let $Q=T[H_1, \dots, H_t]$ be a strong semicomplete composition and $k\geq 4$. The digraph $Q$ has a $k$-kernel if and only if there is a vertex $v\in V(Q)$ such that $\{v\}$ is a $(k-1)$-absorbent set of $Q-(V(H_i)\setminus \{v\})$, where $v\in V(H_i)$ for some $i\in [t]$. 
\end{lem}
\begin{pf}
If there is a vertex $v\in V(Q)$ such that $\{v\}$ is a $(k-1)$-absorbent set of $Q'(=Q-(V(H_i)\setminus \{v\}))$, where $i\in V(H_i)$ for some $i\in [t]$, then $d_Q(u, v)\leq d_{Q'}(u, v)\leq k-1$ for any vertex $u\in V(Q)\setminus V(H_i)$. By Theorem~\ref{thm02}, $u_i$ belongs to a cycle of length three in $T$, which means that $d_Q(u, v)\leq 3$ for each $u\in V(H_i)\setminus \{v\}$. Hence $\{v\}$ is a $(k-1)$-absorbent set of $Q$, and so it is a $k$-kernel of $Q$. The other direction holds by Lemma~\ref{thmc1}.
\end{pf}

The following result holds by the fact that $d_T(u_j, u_i)=d_{Q'}(u, v)$ for any $u\in V(H_j)$ with $j\in [t] \setminus \{i\}$, where $Q'=Q-(V(H_i)\setminus \{v\})$ such that $\{v\}$ is a $k$-absorbent set of $Q'$.

\begin{lem}\label{lem11}
Let $Q=T[H_1, \dots, H_t]$ be a digraph composition and $v\in V(H_i)$ for some $i\in [t]$. The set $\{v\}$ is a $k$-absorbent set of $Q-(V(H_i)\setminus \{v\})$ if and only if $\{u_i\}$ is a $k$-absorbent set of $T$.
\end{lem}

Recall that the problem {\sc $k$-Kernel} for digraphs is NP-complete when $k\in \{2,3\}$. In fact, we can extend this to any fixed integer $k$ which is divisible by 2 or 3.

\begin{lem}\label{lem-kkernel}
Let $k$ be a fixed integer which is divisible by 2 or 3. The problem {\sc $k$-Kernel} for digraphs is NP-complete.
\end{lem}
\begin{pf} Clearly, the problem belongs to NP. Let $D$ be a digraph. We construct a new digraph $D'$ from $D$ as follows: for every vertex $x$ in $D$, substitute $x$ with a path of length $t$: $x_0, x_1, \dots, x_t$, such that all arcs into $x$ in $D$ go into $x_0$ in $D'$ and all arcs out of $x$ in $D$ now go out of $x_t$.

We prove that $D$ contains a kernel if and only if $D'$ contains a $(2t+2)$-kernel, and we complete the proof by the NP-completeness of the problem {\sc Kernel} (i.e. {\sc 2-kernel})  for digraphs.

Let $K$ be a kernel of $D$. We create a $(2t+2)$-kernel in $D'$ by putting $y_t$ in $K'$ for every $y$ in $K$. Let $x_t, y_t \in K'$, clearly, $x, y\in K$. Observe that $d_{D'}(x_t,y_t)\geq 2t+2$ as $d_{D}(x,y)\geq 2$. Let $u=z_i\in V(D')\setminus K'$. If $z_t\in K'$, then $d_{D'}(u, K')= d_{D'}(u, z_t)\leq t< 2t+1$. Otherwise, as $zy\in A(D)$ for some $y\in K$ (by the definition of $K$), we have $d_{D'}(u, K')\leq d_{D'}(z_0, y_t)\leq 2t+1$. The argument above implies that $K'$ is indeed a $(2t+2)$-kernel in $D'$.

Let $K'$ be a $(2t+2)$-kernel in $D'$. We create a kernel in $D$ by putting $y$ in $K$ for every $y_i$ in $K'$. Note that for each $y\in V(D)$, there are at most one $y_i\in V(K')$ by the definition of $K'$. For every pair of vertices $x, y\in K$ (this sentence means that $x_i,y_j\in K'$ for some $i,j$), $d_D(x, y)\geq 2$ (Otherwise, $d_{D'}(x_i, y_j)\leq 2t+1$ for any pair $i, j$ with $i,j\in [t]$, which produces a contradiction). We further claim that $d_D(u, K)=1$ for each $u\in V(D)\setminus K$  (and for each such $u$, $u_0\in V(D)\setminus K'$). Otherwise, $d_D(u,x)\geq 2$ for each $x\in K$, so $d_{D'}(u_0,x_i)\geq 2t+2$, where $x_i\in K'$ corresponds to $x$. Hence $d_{D'}(u_0, K')\geq 2t+2$, which produces a contradiction. The argument above implies that $K$ is indeed a kernel in $D$.

We can also prove that $D$ contains a 3-kernel if and only if $D'$ contains a $(3t+3)$-kernel, by the NP-completeness of the problem {\sc 3-Kernel} for digraphs. The argument is similar, so we omit the details.


\end{pf}

We now study the complexity of {\sc $k$-Kernel} restricted to semicomplete compositions. In particular, we determine the complexity of this problem restricted to strong semicomplete compositions for any integer $k$ with $k\geq 2$.

\begin{thm}\label{thmd2}
Let $k$ be an integer with $k\geq 2$. The following assertions hold:
\begin{description}
\item[(a)] The problem {\sc $k$-Kernel} restricted to strong semicomplete compositions is NP-complete when $k\in \{2,3\}$, and is polynomial-time solvable when $k\geq 4$.
\item[(b)] The problem {\sc $k$-Kernel} restricted to non-strong semicomplete compositions is NP-complete when $k$ is divisible by 2 or 3.
\end{description}
\end{thm}
\begin{pf}

\noindent{\bf Part (a)}
We first consider the case that $k=3$ and use the NP-completeness of the problem {\sc 3-Kernel}. Let $D$ be a digraph and let $Q=\overrightarrow{C}_3[D_1, D_2, D_3]$, where $D_i\cong D$. 
If $D$ contains a 3-kernel, say $K$, then $Q$ contains a 3-kernel, as we just take the set $K$ in $D_1$. If $Q$ contains a 3-kernel, say $K'$, then we note that $K'\subseteq V(D_i)$ for some $i \in \{1,2,3\}$. This argument implies that $K'$ is also a 3-kernel of $D$. Therefore, $D$ contains a 3-kernel if and only if $Q$ does. As $Q$ is a strong semicomplete composition, the result holds for the case $k=3$. The case that $k=2$ can be proved similarly, using a digon instead of $\overrightarrow{C}_3$, so we omit the details.

Now we consider the case that $k\geq 4$. Let $Q=T[H_1, \dots, H_t]$ be a strong semicomplete composition.
By Lemmas~\ref{thmd1} and~\ref{lem11}, $Q$ has a $k$-kernel if and only if $\{u_i\}$ is a $(k-1)$-absorbent set of $T$ for some $i\in [t]$. It suffices to show that we can decide in polynomial time if $\{u_i\}$ is a $(k-1)$-absorbent set of $T$, where $i\in [t]$. Indeed, we just compute the values of $d_{T}(x, u_i)$ for all $x\in V(T)\setminus \{u_i\}$ and check if they are at most $k-1$. This can be done in polynomial time by the Dijkstra's algorithm (see e.g. Theorem~3.3.7 of \cite{Bang-Jensen-Gutin}).

\noindent{\bf Part (b)}  Let $k$ ($k\geq 2$) be an integer which is divisible by 2 or 3. Let $D$ be a digraph and let $Q=T[D_1, \dots, D_t]$, where $D_i\cong D$ and $T$ is a transitive tournament with vertex set $\{u_i\colon\ i\in [t]\}$. Clearly, $Q$ is a non-strong semicomplete composition.
If $D$ has a $k$-kernel, say $K$, then $Q$ contains a $k$-kernel as we just take the set $K$ in $D_t$. If $Q$ has a $k$-kernel, say $K'$, then we have $K\subseteq V(D_t)$, and hence it is also a $k$-kernel of $D_t$, that is, $D$ contains a $k$-kernel. Therefore, $D$ contains a $k$-kernel if and only if $Q$ does. By Lemma~\ref{lem-kkernel}, the conclusion holds.
\end{pf}

\vskip 4mm

\noindent {\bf Conflict of interest} The authors declare that they have no known competing financial interests or personal
relationships that could have appeared to influence the work reported in this paper.

\vskip 1cm

\noindent {\bf Acknowledgement.} We are thankful to Professor Anders Yeo for
discussions on the complexity of {\sc k-Kernel} problem. Yuefang Sun was supported by National Natural Science Foundation of China under Grant No. 12371352, Zhejiang Provincial Natural Science Foundation of China under Grant No.  LY23A010011 and Yongjiang Talent Introduction Programme of Ningbo under Grant No. 2021B-011-G.

\end{document}